# Extended Nonstandard Neutrosophic Logic, Set, and Probability

## based on

## Extended Nonstandard Analysis

Prof. Florentin Smarandache, PhD, Postdoc
University of New Mexico
Mathematics Department
705 Gurley Ave., Gallup, NM 87301, USA

**Abstract**.

We extend for the second time the Nonstandard Analysis by adding the *left monad closed to the right*, and *right monad closed to the left*, while besides the pierced binad (we introduced in 1998) we add now the *unpierced binad* - all these in order to close the newly extended nonstandard space under nonstandard addition, nonstandard subtraction, nonstandard multiplication, nonstandard division, and nonstandard power operations.

Then, we extend the Nonstandard Neutrosophic Logic, Nonstandard Neutrosophic Set, and Nonstandard Probability on this Extended Nonstandard Analysis space - that we prove it is a nonstandard neutrosophic lattice of first type (endowed with a nonstandard neutrosophic partial order) as well as a nonstandard neutrosophic lattice of second type (as algebraic structure, endowed with two binary neutrosophic laws, $inf_N$ and $sup_N$).

Many theorems, new terms introduced, and examples of nonstandard neutrosophic operations are given.

## 1. Short Introduction

In order to more accurately situate and fit the neutrosophic logic into the framework of extended nonstandard analysis, we present the nonstandard neutrosophic inequalities, nonstandard neutrosophic equality, nonstandard neutrosophic infimum and supremum, nonstandard neutrosophic intervals, including the cases when the neutrosophic logic standard and nonstandard components *T, I, F* get values outside of the classical unit interval *[0, 1],* and a brief evolution of neutrosophic operators.

## 2. Theoretical Reason for the Nonstandard Form of Neutrosophic Logic

The only reason I have added the nonstandard form to neutrosophic logic (and similarly to neutrosophic set and probability) was in order to make a distinction between *Relative Truth* (which is truth in some Worlds, according to Leibniz) and *Absolute Truth* (which is truth in all possible Words, according to Leibniz as well) that occur in philosophy.

Another possible reason may be when the neutrosophic degrees of truth, indeterminacy, or falsehood are infinitesimally determined, for example a value infinitesimally bigger than *0.8 (*or *0.8⁺)*, or infinitesimally smaller than *0.8 (*or *⁻0.8)*. But these can easily be overcome by roughly using interval neutrosophic values, for example (0.80, 0.81) and (0.79, 0.80) respectively.

### 3. Why the Sum of Neutrosophic Components is up to 3

I was more prudent when I presented the sum of single valued standard neutrosophic components, saying:

*Let T, I, F be single valued numbers, T, I, F ∈ [0, 1], such that $0 \leq T + I + F \leq 3$.* (1)

A friend alerted me: *"If T, I, F are numbers in [0, 1], of course their sum is between 0 and 3."* "Yes, I responded, I afford this tautology, because if I did not mention that the sum is up to *3*, readers would take for granted that the sum $T + I + F$ is bounded by *1*, since that is in all previous logics and in probability!"

Their sum is up to 3 since they are considered independent of each other.

### 4. Neutrosophic Components Outside the Unit Interval [0, 1]

Thinking out of box, inspired from the real world, was the first intent, i.e. allowing neutrosophic components (truth / indeterminacy / falsehood) values be outside of the classical (standard) unit real interval *[0, 1]* used in all previous (Boolean, multi-valued etc.) logics if needed in applications, so neutrosophic component values *< 0* and *> 1* had to occurs due to the Relative / Absolute stuff, with:

$$^-0 <_N 0 \quad \text{and} \quad 1^+ >_N 1.$$  (2)

Later on, in 2007, I found plenty of cases and real applications in Standard Neutrosophic Logic and Set (therefore, not using the Nonstandard Neutrosophic Logic and Set), and it was thus possible the extension of the neutrosophic set *to Neutrosophic Overset (when some neutrosophic component is > 1), and to Neutrosophic Underset (when some neutrosophic component is < 0), and to Neutrosophic Offset (when some neutrosophic components are off the interval [0, 1], i.e. some neutrosophic component > 1 and some neutrosophic component < 0). Then, similar extensions to respectively Neutrosophic Over/Under/Off Logic, Measure, Probability, Statistics etc.* [8, 17, 18, 19], extending the unit interval [0, 1] to

*[Ψ, Ω], with $\Psi \leq 0 < 1 \leq \Omega$,* (3)

where *Ψ, Ω* are standard real numbers.

### 5. Refined Neutrosophic Logic, Set, and Probability

I wanted to get the neutrosophic logic as general as possible [6], extending all previous logics (Boolean, fuzzy, intuitionistic fuzzy logic, intuitionistic logic, paraconsistent logic, dialethism), and to have it able to deal with all kind of logical propositions (including paradoxes, nonsensical propositions, etc.).

That's why in 2013 I extended the Neutrosophic Logic to *Refined Neutrosophic Logic* [ from generalizations of *2-valued Boolean logic* to fuzzy logic, also from the *Kleene's and Lukasiewicz's* and *Bochvar's 3-symbol valued logics* or *Belnap's 4-symbol valued logic* to the most general *n-symbol* or *n-numerical* valued refined neutrosophic logic, for any integer $n \geq 1$ ], the largest ever so far, when some or all neutrosophic components *T, I, F* were respectively split/refined into neutrosophic subcomponents: $T_1, T_2, \ldots; I_1, I_2, \ldots; F_1, F_2, \ldots$ which were deduced from our everyday life [3].

6. **From Paradoxism movement to Neutrosophy branch of philosophy and then to Neutrosophic Logic**

   I started first from *Paradoxism* (that I founded in 1980's as a movement based on antitheses, antinomies, paradoxes, contradictions in literature, arts, and sciences), then I introduced the *Neutrosophy* (as generalization of Dialectics of Hegel and Marx, which is actually the ancient YinYang Chinese philosophy), neutrosophy is a branch of philosophy studying the dynamics of triads, inspired from our everyday life, triads that have the form:
   
       *<A>*, its opposite *<antiA>*, and their neutrals *<neutA>*,                                    (4)
   
   where *<A>* is any item or entity [4].
   (Of course, we take into consideration only those triads that make sense in our real and scientific world.)

   The Relative Truth neutrosophic value was marked as *1*, while the Absolute Truth neutrosophic value was marked as $1^+$ (a tinny bigger than the Relative Truth's value): $1^+ >_N 1$, where $>_N$ is a neutrosophic inequality, meaning $1^+$ is neutrosophically bigger than *1*.
   Similarly for Relative Falsehood / Indeterminacy (which falsehood / indeterminacy in some Worlds), and Absolute Falsehood / Indeterminacy (which is falsehood / indeterminacy in all possible worlds).

7. **Introduction to Nonstandard Analysis** [15, 16]

   An <u>*infinitesimal*</u> [or infinitesimal number] ($\varepsilon$) is a number $\varepsilon$ such that $|\varepsilon| < 1/n$, for any non-null positive integer *n*. An infinitesimal is close to zero, and so small that it cannot be measured.

The infinitesimal is a number smaller, in absolute value, than anything positive nonzero.

Infinitesimals are used in calculus.

    An <u>*infinite*</u> [or infinite number] *( $\omega$ )* is a number greater than anything:

$1 + 1 + 1 + \ldots + 1$ (for any finite number terms)                                                                (5)

The infinites are reciprocals of infinitesimals.

The set of *hyperreals* (or *non-standard reals*), denoted as $R^*$, is the extension of set of the real numbers, denoted as $R$, and it comprises the infinitesimals and the infinites, that may be represented on the *hyperreal number line*

$$1/\varepsilon = \omega/1. \qquad (6)$$

The set of hyperreals satisfies the *transfer principle*, which states that the statements of first order in $R$ are valid in $R^*$ as well.

A *monad* (*halo*) of an element $a \in R^*$, denoted by $\mu(a)$, is a subset of numbers infinitesimally close to $a$.

### 8. First Extension of Nonstandard Analysis

Let's denote by $R_+^*$ the set of positive nonzero hyperreal numbers.

We consider the left monad and right monad, and the (*pierced*) binad that we have introduced as extension in 1998 [5]:

**Left Monad** { that we denote, for simplicity, by $(^-a)$ or only $^-a$ } is defined as:

$$\mu(^-a) = (^-a) = {}^-a = \bar{a} = \{a - x,\ x \in R_+^* \mid x \text{ is infinitesimal}\}. \qquad (7)$$

**Right Monad** { that we denote, for simplicity, by $(a^+)$ or only by $a^+$ } is defined as:

$$\mu(a^+) = (a^+) = a^+ = \overset{+}{a} = \{a + x,\ x \in R_+^* \mid x \text{ is infinitesimal}\}. \qquad (8)$$

**Pierced Binad** { that we denote, for simplicity, by $(^-a^+)$ or only $^-a^+$ } is defined as:

$$\mu(^-a^+) = (^-a^+) = {}^-a^+ = \overset{-+}{a} =$$

$$= \{a - x,\ x \in R_+^* \mid x \text{ is infinitesimal}\} \cup \{a + x,\ x \in R_+^* \mid x \text{ is infinitesimal}\}$$

$$= \{a \pm x,\ x \in R_+^* \mid x \text{ is infinitesimal}\}. \qquad (9)$$

The left monad, right monad, and the pierced binad are subsets of $R^*$.

### 9. Second Extension of Nonstandard Analysis

For necessity of doing calculations that will be used in nonstandard neutrosophic logic in order to calculate the nonstandard neutrosophic logic operators (conjunction, disjunction, negation, implication, equivalence) and in order to have the Nonstandard Real MoBiNad Set closed under arithmetic operations, we extend now for the time: the left monad to the Left Monad Closed to the Right, the right monad to the Right Monad Closed to the Left; and the Pierced Binad to the Unpierced Binad, defined as follows:

**Left Monad Closed to the Right**

$$\mu\!\left(\overset{-0}{a}\right) = \left(\overset{-0}{a}\right) = \overset{-0}{a} = \{a - x \mid x = 0,\ \text{or}\ x \in R_+^* \text{ and } x \text{ is infinitesimal}\} = \mu(^-a) \cup \{a\} = (^-a) \cup \{a\}$$

$$= {}^-a \cup \{a\}. \tag{10}$$

**Right Monad Closed to the Left**

$$\mu\left(\overset{0+}{a}\right) = \left(\overset{0+}{a}\right) = \overset{0+}{a} = \{a + x \mid x = 0, \text{ or } x \in R_+^* \text{ and } x \text{ is infinitesimal}\} = \mu(a^+) \cup \{a\} = (a^+) \cup \{a\}$$

$$= a^+ \cup \{a\}. \tag{11}$$

**Pierced Binad**

$$\mu\left(\overset{-0+}{a}\right) = \left(\overset{-0+}{a}\right) = \overset{-0+}{a} = \{a - x \mid x \in R_+^* \text{ and } x \text{ is infinitesimal}\} \cup \{a + x \mid x \in R_+^* \text{ and } x \text{ is infinitesimal}\} \cup \{a\} =$$

$$= \{a \pm x \mid x = 0, \text{ or } x \in R_+^* \text{ and } x \text{ is infinitesimal}\}$$

$$= \mu({}^-a^+) \cup \{a\} = ({}^-a^+) \cup \{a\} = {}^-a^+ \cup \{a\} \tag{12}$$

The element $\{a\}$ has been included into the left monad, right monad, and pierced binad respectively.

### 10. Nonstandard Neutrosophic Function

In order to be able to define equalities and inequalities in the sets of monads, and in the sets of binads, we construct a *nonstandard neutrosophic function* that approximates the monads and binads to tiny open (or half open and half closed respectively) standard real intervals as below. It is called 'neutrosophic' since it deals with indeterminacy: unclear, vague monads and binads, and the function approximates them with some tiny real subsets.

Taking an arbitrary infinitesimal
$\varepsilon_1 > 0$, and writing ${}^-a = a - \varepsilon_1$, $a^+ = a + \varepsilon_1$, and ${}^-a^+ = a \pm \varepsilon_1$, (13)

or taking an arbitrary infinitesimal $\varepsilon_2 \geq 0$, and writing

$$\overset{-0}{a} = (a - \varepsilon_2, a], \overset{0+}{a} = [a, a + \varepsilon_2), \overset{-0+}{a} = (a - \varepsilon_2, a + \varepsilon_2) \tag{14}$$

we meant actually picking up a representative from each class of the monads and of the binads respectively.

Representations of the monads and binads by intervals is not quite accurate from a classical point of view, but it is an approximation that helps in finding a partial order and computing nonstandard arithmetic operations on the elements of the nonstandard set $NR_{MB}$.

Let $\varepsilon$ be a generic positive infinitesimal, while $a$ be a generic standard real number.

Let $P(R)$ be the power set of the real number set $R$.

$$\mu_N: NR_{MB} \rightarrow P(R) \tag{15}$$

For any $a \in R$, the set of real numbers, one has:

$$\mu_N(({}^-a)) =_N (a - \varepsilon, a), \tag{16}$$

$$\mu_N((a^+)) =_N (a, a + \varepsilon), \tag{17}$$

$$\mu_N(\ ({}^-a^+)\ ) =_N (a-\varepsilon, a) \cup (a, a+\varepsilon),\qquad(18)$$

$$\mu_N\left(\begin{pmatrix}{}^{-0}a\end{pmatrix}\right) =_N (a-\varepsilon, a],$$

$$\mu_N\left(\begin{pmatrix}{}^{0+}a\end{pmatrix}\right) =_N [a, a+\varepsilon),\qquad(19\text{-}20\text{-}21)$$

$$\mu_N\left(\begin{pmatrix}{}^{-0+}a\end{pmatrix}\right) =_N (a-\varepsilon, a+\varepsilon),$$

$$\mu_N\left(\begin{pmatrix}{}^{0}a\end{pmatrix}\right) =_N \mu_N(a) =_N a = [a,a],\qquad(22)$$

in order to set it as real interval too.

## 11. General Notations for Monads and Binads

Let $a \in R$ be a standard real number. We use the following general notation for monads and binads:

$${}^m a \in \{{}^-a, {}^{-0}a, {}^+a, {}^{0+}a, {}^{-+}a, {}^{-0+}a,\ {}^0 a\}\text{ and by convention }{}^0 a = a;\qquad(23)$$

$$\text{or } m \in \{\ ,\ ^-,\ ^{-0},\ ^+,\ ^{+0},\ ^{-+},\ ^{-0+}\} = \{{}^0,\ ^-,\ ^{-0},\ ^+,\ ^{+0},\ ^{-+},\ ^{-0+}\};\qquad(24)$$

therefore "$m$" above a standard real number "$a$" may mean anything: a standard real number ($^0$, or nothing above), a left monad ($^-$), a left monad closed to the right ($^{-0}$), a right monad ($^+$), a right monad closed to the left ($^{0+}$), a pierced binad ($^{-+}$), or a unpierced binad ($^{-0+}$) respectively.

The notations of monad's and binad's diacritics above (not laterally) the number $a$ as

$${}^-a, {}^{-0}a, {}^+a, {}^{0+}a, {}^{-+}a, {}^{-0+}a\qquad(25)$$

are the best, since they also are designed to avoid confusion for the case when the real number $a$ is negative.

For example, if $a = {}^-2$, then the corresponding monads and binads are respectively represented as: ${}^-\!\!{}^-2, {}^{-0}\!\!{}^-2, {}^+\!\!{}^-2, {}^{0+}\!\!{}^-2, {}^{-+}\!\!{}^-2, {}^{-0+}\!\!{}^-2$. (26)

## 12. Classical and Neutrosophic Notations

Classical notations on the set of real numbers:
$<, \leq, >, \geq, \wedge, \vee, \rightarrow, \leftrightarrow, \cap, \cup, \subset, \supset, \subseteq, \supseteq, =, \in,$
$+, -, \times, \div, \wedge, *$ (27)
Operations with real subsets: $\circledast$ (28)
Neutrosophic notations on nonstandard sets (that involve indeterminacies, approximations, vague boundaries):

$$<_N, \leq_N, >_N, \geq_N, \wedge_N, \vee_N, \rightarrow_N, \leftrightarrow_N, \cap_N, \cup_N, \subset_N, \supset_N, \subseteq_N, \supseteq_N, =_N, \in_N$$
$$+_N, -_N, \times_N, \div_N, \wedge_N, *_N \tag{29}$$

### 13. Neutrosophic Strict Inequalities

We recall the neutrosophic strict inequality which is needed for the inequalities of nonstandard numbers.

Let α, β be elements in a partially ordered set *M*.

We have defined the *neutrosophic strict inequality*

$$\alpha >_N \beta \tag{30}$$

and read as

"*α is neutrosophically greater than β*"

if

*α in general is greater than β,*

or *α is approximately greater than β,*

or *subject to some indeterminacy* (unknown or unclear ordering relationship between α and β) *or subject to some contradiction (*situation when α is smaller than or equal to β*) α is greater than β.*

It means that in most of the cases, on the set *M*, *α is greater than β*.

And similarly for the opposite neutrosophic strict inequality $\alpha <_N \beta$. (31)

### 14. Neutrosophic Equality

We have defined the *neutrosophic inequality*

$$\alpha =_N \beta \tag{32}$$

and read as

"*α is neutrosophically equal to β*"

if

*α in general is equal to β,*

or *α is approximately equal to β,*

or *subject to some indeterminacy* (unknown or unclear ordering relationship between α and β) *or subject to some contradiction (*situation when α is not equal to β*) α is equal to β.*

It means that in most of the cases, on the set *M*, *α is equal to β*.

### 15. Neutrosophic (Non-Strict) Inequalities

Combining the neutrosophic strict inequalities with neutrosophic equality, we get the $\geq_N$ and $\leq_N$ neutrosophic inequalities.

Let $\alpha, \beta$ be elements in a partially ordered set $M$.

The *neutrosophic (non-strict) inequality*

$$\alpha \geq_N \beta \tag{33}$$

and read as

*"α is neutrosophically greater than or equal to β"*

if

*α in general is greater than or equal to β,*

or *α is approximately greater than or equal to β,*

or *subject to some indeterminacy* (unknown or unclear ordering relationship between α and β) *or subject to some contradiction (*situation when α is smaller than β*) α is greater than or equal to β.*

It means that in most of the cases, on the set $M$, $\alpha$ is greater than or equal to $\beta$.

And similarly for the opposite neutrosophic (non-strict) inequality $\alpha \leq_N \beta$. (34)

### 16. Neutrosophically Ordered Set

Let $M$ be a set. $(M, <_N)$ is called a neutrosophically ordered set if:

$$\forall\ \alpha, \beta \in M,\text{ one has: either } \alpha <_N \beta,\text{ or } \alpha =_N \beta,\text{ or } \alpha >_N \beta. \tag{35}$$

### 17. Neutrosophic Infimum and Neutrosophic Supremum

As an extension of the classical infimum and classical supremum, and using the neutrosophic inequalities and neutrosophic equalities, we define the neutrosophic infimum ( denoted as $inf_N$ ) and the neutrosophic supremum ( denoted as $sup_N$ ).

*Neutrosophic Infimum.*
Let $(S, <_N)$ be a set that is neutrosophically partially ordered, and $M$ a subset of $S$. The neutrosophic infimum of $M$, denoted as $inf_N(M)$ is the neutrosophically greatest element in $S$ that is neutrosophically less than or equal to all elements of $M$:

*Neutrosophic Supremum.*
Let $(S, <_N)$ be a set that is neutrosophically partially ordered, and $M$ a subset of $S$. The neutrosophic supremum of $M$, denoted as $sup_N(M)$ is the neutrosophically smallest element in $S$ that is neutrosophically greater than or equal to all elements of $M$.

## 18. Definition of Nonstandard Real MoBiNad Set

Let $\mathbb{R}$ be the set of standard real numbers, $\mathbb{R}^*$ the set of hyper-reals (or non-standard reals) which consists of infinitesimals and infinites.

The Nonstandard Real MoBiNad Set is now defined, for the first time, as follows:

$$NR_{MB} =_N \left\{ \begin{array}{l} \varepsilon, \omega, a, (^-a), (^-a^0), (a^+), (^0a^+), (^-a^+), (^-a^{0+}) \text{ |where } \varepsilon \text{ are infinitesimals,} \\ \text{with } \varepsilon \in \mathbb{R}^*; \omega \text{ are infinites, with } \omega \in \mathbb{R}^*; \text{ and } a \text{ are real numbers, with } a \in \mathbb{R} \end{array} \right\}. \quad (36)$$

Therefore:

$$NR_{MB} =_N \mathbb{R}^* \cup \mathbb{R} \cup \mu(^-\mathbb{R}) \cup \mu(^-\mathbb{R}^0) \cup \mu(\mathbb{R}^+) \cup \mu(^0\mathbb{R}^+) \cup \mu(^-\mathbb{R}^+) \cup \mu(^-\mathbb{R}^{0+}), \quad (37)$$

where $\mu(^-\mathbb{R})$ is the set of all real left monads,

$\mu(^-\mathbb{R}^0)$ is the set of all real left monads closed to the right,

$\mu(\mathbb{R}^+)$ is the set of all real right monads,

$\mu(^0\mathbb{R}^+)$ is the set of all real right monads closed to the left,

$\mu(^-\mathbb{R}^+)$ is the set of all real pierced binads,

and $\mu(^-\mathbb{R}^{0+})$ is the set of all real unpierced binads.

Also, $NR_{MB} =_N \left\{ \varepsilon, \omega, \overset{m}{a} \text{ |where } \varepsilon, \omega \in \mathbb{R}^*; a \in \mathbb{R}; \text{ and } m \in \{^-, ^{-0}, ^+, ^{+0}, ^{-+}, ^{-0+}\} \right\}. \quad (38)$

## 19. Etymology of MoBiNad

MoBiNad comes from **mo**nad + **bi**nad, introduced now for the first time.

## 20. Definition of Nonstandard Complex MoBiNad Set

The **Nonstandard Complex MoBiNad Set**, introduced here for the first time, is defined as:

$$NC_{MB} =_N \{\alpha + \beta i | \text{ where } i = \sqrt{-1}; \alpha, \beta \in NR_{MB}\}. \quad (39)$$

## 21. Definition of Nonstandard Neutrosophic Real MoBiNad Set

The **Nonstandard Neutrosophic Real MoBiNad Set**, introduced now for the first time, is defined as:

$$NNR_{MB} =_N \{\alpha + \beta I | \text{ where } I = \text{literal indeterminacy}, I^2 = I; \alpha, \beta \in NR_{MB}\}. \quad (40)$$

## 22. Definition of Nonstandard Neutrosophic Complex MoBiNad Set

The **Nonstandard Neutrosophic Complex MoBiNad Set**, introduced now for the first time, is defined as:

$NNC_{MB} =_N \{\alpha + \beta I |$ where $I =$ literal indeterminacy, $I^2 = I$; $\alpha, \beta \in NC_{MB}\}$.  (41)

### 23. Properties of the Nonstandard Neutrosophic Real Mobinad Set

Since in nonstandard neutrosophic logic we use only the nonstandard neutrosophic real mobinad set, we study some properties of it.

**Theorem 1**

The nonstandard real mobinad set $(NR_{MB}, \leq_N)$, endowed with the nonstandard neutrosophic inequality is a lattice of first type [as partially ordered set (*poset*)].

*Proof*

The set $NR_{MB}$ is partially ordered, because [except the two-element subsets of the form $\{a, {}^{-+}a\}$, and $\{a, {}^{-0+}a\}$, with $a \in \mathbb{R}$, beetwen which there is no order] all other elements are ordered:

If $a < b$, where $a, b \in \mathbb{R}$, then: ${}^{m_1}a <_N {}^{m_2}b$, for any monads or binads

$m_1, m_2 \in_N \{{}^-, {}^{-0}, {}^+, {}^{0+}, {}^{-+}, {}^{-0+}\}$.  (42)

If $a = b$, one has:

${}^-a <_N a$,  (43)

$a^- <_N a^+$,  (44)

$a <_N a^+$,  (45)

${}^-a \leq_N {}^-a^+$,  (46)

${}^-a^+ \leq_N a^+$,  (47)

and there is no neutrosophic ordering relationship between $a$ and ${}^-a^+$,

nor between $a$ and ${}^{-0+}a$ (that is why $\leq_N$ on $NR_{MB}$ is a partial ordering set).  (48)

If $a > b$, then: ${}^{m_1}a >_N {}^{m_2}b$, for any monads or binads $m_1, m_2$.  (49)

Any two-element set $\{\alpha, \beta\} \subset_N NR_{MB}$ has a **neutrosophic nonstandard infimum** (**meet**, or **greatest lower bound**) that we denote by $\inf_N$, and a **neutrosophic nonstandard supremum** (**joint**, or **least upper bound**) that we denote by $\sup_N$, where both

$\inf_N\{\alpha, \beta\}$ and $\sup_N\{\alpha, \beta\} \in NR_{MB}$.  (50)

For the non-ordered elements $a$ and $^-a^+$:

$$\inf_N\{a,^- a^+\} =_N {^-a} \in_N NR_{MB}, \tag{51}$$

$$\sup_N\{a,^- a^+\} =_N a^+ \in_N NR_{MB}. \tag{52}$$

And similarly for non-ordered elements $a$ and $^-a^{0\,+}$:

$$\inf_N\{a,^- a^{0\,+}\} =_N {^-a} \in_N NR_{MB}, \tag{53}$$

$$\sup_N\{a,^- a^{0\,+}\} =_N a^+ \in_N NR_{MB}. \tag{54}$$

Dealing with monads and binads which neutrosophically are real subsets with indeterminate borders, and similarly $a = [a, a]$ can be treated as a subset, we may compute $\inf_N$ and $\sup_N$ of each of them.

$$\inf_N(^-a) =_N {^-a} \text{ and } \sup_N(^-a) =_N {^-a}; \tag{55}$$

$$\inf_N(a^+) =_N a^+ \text{ and } \sup_N(a^+) =_N a^+; \tag{56}$$

$$\inf_N(^-a^+) =_N {^-a} \text{ and } \sup_N(^-a^+) =_N a^+; \tag{57}$$

$$\inf_N(^-a^{0\,+}) =_N {^-a} \text{ and } \sup_N(^-a^{0\,+}) =_N a^+. \tag{58}$$

Also, $\inf_N(a) =_N a$ and $\sup_N(a) =_N a$. \tag{59}

If $a < b$, then $\overset{m_1}{a} <_N \overset{m_2}{b}$, whence $\inf_N\left\{\overset{m_1}{a}, \overset{m_2}{b}\right\} =_N \inf_N\left(\overset{m_1}{a}\right)$ and $\sup_N\left\{\overset{m_1}{a}, \overset{m_2}{b}\right\} =_N \sup_N\left(\overset{m_2}{b}\right)$,

$$\tag{60}$$

which are computed as above.

Similarly, if $a > b$, with $\overset{m_1}{a} >_N \overset{m_2}{b}$. \tag{61}

If $a = b$, then:

$\inf_N\left\{\overset{m_1}{a}, \overset{m_2}{a}\right\} =_N$ the neutrosophically smallest ($<_N$) element among $\inf_N\left\{\overset{m_1}{a}\right\}$ and $\inf_N\left\{\overset{m_2}{a}\right\}$.

$$\tag{62}$$

While $\sup_N\left\{\overset{m_1}{a}, \overset{m_2}{a}\right\} =_N$ the neutrosophically greatest ($>_N$) element among $\sup_N\left\{\overset{m_1}{a}\right\}$ and $\sup_N\left\{\overset{m_2}{a}\right\}$. \tag{63}

Examples:

$$\inf_N(^-a, a^+) =_N {^-a} \text{ and } \sup_N(^-a, a^+) =_N a^+; \tag{64}$$

$$\inf{}_N({}^-a, {}^-a^+) =_N {}^-a \text{ and } \sup{}_N({}^-a, {}^-a^+) =_N {}_N a^+; \qquad (65)$$

$$\inf{}_N({}^-a^+, a^+) =_N {}^-a \text{ and } \sup{}_N({}^-a^+, a^+) =_N a^+. \qquad (66)$$

Therefore, $(NR_{MB}, \leq_N)$ is a nonstandard real mobinad lattice of first type (as partially ordered set).

**Consequence**

If we remove all pierced and unpierced binads from $NR_{MB}$ and we denote the new set by

$NR_M = \{\varepsilon, \omega, a, {}^-a, {}^-a^0, a^+, {}^0a^+, \text{where } \varepsilon \text{ are infinitesimals}, \omega \text{ are infinites, and } a \in \mathbb{R}\}$ we obtain a totally neutrosophically ordered set.

**Theorem 2**

Any finite non-empty subset $L$ of $(NR_{MB}, \leq_N)$ is also a sublattice of first type.

*Proof*

It is a consequence of any classical lattice of first order (as partially ordered set).

**Theorem 3**

$(NR_{MB}, \leq_N)$ is not bounded neither to the left nor to the right, since it does not have a **minimum** (bottom, or least element), nor a **maximum** (top, or greatest element).

**Theorem 4**

$(NR_{MB}, \inf_N, \sup_N)$, where $\inf_N$ and $\sup_N$ are two binary operations, dual to each other, defined before, is a **lattice of second type** (as an algebraic structure).

*Proof*

We have to show that the two laws $\inf_N$ and $\sup_N$ are commutative, associative, and verify the absorption laws.

Let $\alpha, \beta, \gamma \in NR_{MB}$ be two arbitrary elements.

**Commutativity Laws**

$$\text{i) } \inf{}_N\{\alpha, \beta\} =_N \inf{}_N\{\beta, \alpha\}. \qquad (67)$$

$$\text{ii) } \sup{}_N\{\alpha, \beta\} =_N \sup{}_N\{\beta, \alpha\}. \qquad (68)$$

Their proofs are straightforward.

**Associativity Laws**

i) $\quad \inf_N\{\alpha, \inf_N\{\beta, \gamma\}\} =_N \inf_N\{\inf_N\{\alpha, \beta\}, \gamma\}.$ (69)

*Proof*

$$\inf_N\{\alpha, \inf_N\{\beta, \gamma\}\} =_N \inf_N\{\alpha, \beta, \gamma\}, \tag{70}$$

and

$$\inf_N\{\inf_N\{\alpha, \beta\}, \gamma\} =_N \inf_N\{\alpha, \beta, \gamma\}, \tag{71}$$

where we have extended the binary operation $\inf_N$ to a trinary operation $\inf_N$.

ii) $\sup_N\{\alpha, \sup_N\{\beta, \gamma\}\} =_N \sup_N\{\sup_N\{\alpha, \beta\}, \gamma\}$ (72)

*Proof*

$$\sup_N\{\alpha, \sup_N\{\beta, \gamma\}\} =_N \sup_N\{\alpha, \beta, \gamma\}, \tag{73}$$

and

$$\sup_N\{\sup_N\{\alpha, \beta\}, \gamma\} =_N \sup_N\{\alpha, \beta, \gamma\}, \tag{74}$$

where similarly we have extended the binary operation $\sup_N$ to a trinary operation $\sup_N$.

**Absorption Laws (as peculiar axioms to the theory of lattice)**

i) We need to prove that $\inf_N\{\alpha, \sup_N\{\alpha, \beta\}\} =_N \alpha.$ (75)

Let $\alpha \leq_N \beta$, then $\inf_N\{\alpha, \sup_N\{\alpha, \beta\}\} =_N \inf_N\{\alpha, \beta\} =_N \alpha.$ (76)

Let $\alpha >_N \beta$, then $\inf_N\{\alpha, \sup_N\{\alpha, \beta\}\} =_N \inf_N\{\alpha, \alpha\} =_N \alpha.$ (77)

ii) Now, we need to prove that $\sup_N\{\alpha, \inf_N\{\alpha, \beta\}\} =_N \alpha.$ (78)

Let $\alpha \leq_N \beta$, then $\sup_N\{\alpha, \inf_N\{\alpha, \beta\}\} =_N \sup_N\{\alpha, \alpha\} =_N \alpha.$ (79)

Let $\alpha >_N \beta$, then $\sup_N\{\alpha, \inf_N\{\alpha, \beta\}\} =_N \sup_N\{\alpha, \beta\} =_N \alpha.$ (80)

**Consequence**

The binary operations $\inf_N$ and $\sup_N$ also satisfy the idempotent laws:

$$\inf_N\{\alpha, \alpha\} =_N \alpha, \tag{81}$$

$$\sup_N\{\alpha, \alpha\} =_N \alpha. \tag{82}$$

*Proof*

The axioms of idempotency follow directly from the axioms of absorption proved above.

Thus, we have proved that *(NR_MB, inf_N, sup_N)* is a lattice of second type (as algebraic structure).

## 24. Definition of General Nonstandard Real MoBiNad Interval

Let $a, b \in \mathbb{R}$, with $-\infty < a \leq b < \infty$, (83)

$]^-a, b^+[_{MB} = \{x \in NR_{MB}, {}^-a \leq_N x \leq_N b^+\}$. (84)

As particular edge cases:

$]^-a, a^+[_{MB} =_N \{{}^-a, a, {}^-a^+, a^+\}$, a discrete nonstandard real set of cardinality 4. (85)

$]^-a, {}^-a[_{MB} =_N \{{}^-a\}$; (86)

$]a^+, a^+[_{MB} =_N \{a^+\}$; (87)

$]a, a^+[_{MB} =_N \{a, a^+\}$; (88)

$]^-a, a[_{MB} =_N \{{}^-a, a\}$; (89)

$]^-a, {}^-a^+[_{MB} =_N \{{}^-a, {}^-a^+, a^+\}$, where $a \notin ]^-a, {}^-a^+[_{MB}$ since $a \not\leq_N {}^-a^+$ (there is no relation of order between $a$ and ${}^-a^+$); (90)

$]^-a^+, a^+[_{MB} =_N \{{}^-a^+, a^+\}$. (91)

**Theorem 5**

$(]^-a, b^+[, \leq_N)$ is a nonstandard real mobinad sublattice of first type (poset). (92)

**Theorem 6**

$(]^-a, b^+[, \inf_N, \sup_N, {}^-a, b^+)$ is a nonstandard *bounded* real mobinad sublattice of second type (as algebraic structure). (93)

*Proofs*

$]^-a, b^+[_{MB}$ as a nonstandard subset of $NR_{MB}$ is also a poset, and for any two-element subset $\{\alpha, \beta\} \subset_N ]^-0, 1^+[_{MB}$ (94)

one obviously has the triple neutrosophic nonstandard inequality:

${}^-a \leq_N \inf_N\{\alpha, \beta\} \leq_N \sup_N\{\alpha, \beta\} \leq_N b^+$, (95)

whence $(]^-a, b^+[_{MB} \leq_N)$ is a nonstandard real mobinad sublattice of first type (poset), or sublattice of $NR_{MB}$.

Further on, $]^-a, b^+[$, endowed with two binary operations $\inf_N$ and $\sup_N$, is also a sublattice of the lattice $NR_{MB}$, since the lattice axioms (Commutative Laws, Associative Laws, Absortion Laws, and Idempotent Laws) are clearly verified on $]^-a, b^+[$.

The nonstandard neutrosophic modinad **Identity Join Element** (Bottom) is $^-a$, and the nonstandard neutrosophic modinad **Identity Meet Element** (Top) is $b^+$,

or $\inf_N\,]^-a, b^+[\, =_N\, ^-a$ and $\sup_N\,]^-a, b^+[\, =_N\, b^+$. (96)

The sublattice **Identity Laws** are verified below.

Let $\alpha \in_N\,]^-a, b^+[$, whence $^-a \leq_N \alpha \leq_N b^+$. (97)

Then:

$\inf_N\{\alpha, b^+\} =_N \alpha$, and $\sup_N\{\alpha, ^-a\} =_N \alpha$. (98)

### 25. Definition of Nonstandard Real MoBiNad Unit Interval

$]^-0, 1^+[_{MB} =_N \{x \in NR_{MB}, ^-0 \leq_N x \leq_N 1^+\}$ (99)

$=_N \left\{ \overset{-}{\varepsilon}, \overset{-0}{a}, \overset{+}{a}, \overset{0+}{a}, \overset{-+}{a}, \overset{-0+}{a}, a \,\middle|\, \begin{array}{l} \text{where } \varepsilon \text{ are infinitesimals,} \\ \varepsilon \in \mathbb{R}^*, \text{with } \varepsilon > 0, \text{and } a \in [0, 1] \end{array} \right\}$. (100)

This is an extension of the previous definition (1998) of nonstandard unit interval

$]^-0, 1^+[\, =_N (^-0) \cup [0, 1] \cup (1^+)$ (101)

Associated to the first published definitions of neutrosophic set, logic, and probability was used.

One has: $]^-0, 1^+[\, \subset_N \,]^-0, 1^+[_{MB}$, (102)

where the index $_{MB}$ means: all monads and binads included in $]^-0, 1^+[$, for example:

$(^-0.2), (^-0.3^0), (0.5^+), (^-0.7^+), (^-0.8^{0+})$ etc. (103)

{or, using the top diacritics notation, respectively: $\overset{-}{0.2}, \overset{-0}{0.3}, \overset{+}{0.5}, \overset{-+}{0.7}, \overset{-0+}{0.8}$ etc.}. (104)

**Theorem 7**

The Nonstandard Real MoBiNad Unit Interval $]^-0, 1^+[_{MB}$ is a partially ordered set (poset) with respect to $\leq_N$, and any of its two elements have an $\inf_N$ and $\sup_N$ whence $]^-0, 1^+[_{MB}$ is a nonstandard neutrosophic lattice of first type (as poset).

**Theorem 8**

The Nonstandard Real MoBiNad Unit Interval $]^-0, 1^+[_{MB}$, endowed with two binary operations $\inf_N$ and $\sup_N$, is also a nonstandard neutrosophic lattice of second type (as an algebraic structure).

*Proofs*

Replace $a = 0$ and $b = 1$ into the general nonstandard real mobinad interval $]^-a, b^+[$.

### 26. Definition of Extended General Neutrosophic Logic

We extend and present in a clearer way our 1995 definition (published in 1998) of neutrosophic logic.

Let $\mathcal{U}$ be a universe of discourse of propositions, and $P \in \mathcal{U}$ a generic proposition.

A **General Neutrosophic Logic** is a multivalued logic in which each proposition $P$ has a degree of truth $(T)$, a degree of indeterminacy $(I)$, and a degree of falsehood $(F)$, where $T, I, F$ are standard or nonstandard real mobinad subsets of the nonstandard real mobinat unit interval $]^-0, 1^+[_{MB}$,

with $T, I, F \subseteq_N ]^-0, 1^+[_{MB}$, $\hspace{2cm}$ (105)

where $^-0 \leq_N \inf_N T + \inf_N I + \inf_N F \leq_N \sup_N T + \sup_N I + \sup_N F \leq 3^+$. $\hspace{1cm}$ (106)

### 27. Definition of Standard Neutrosophic Logic

If in the above definition of general neutrosophic logic all neutrosophic components, $T, I, F$, are standard real subsets, included in or equal to the standard real unit interval, $T, I, F \subseteq [0, 1]$, where $0 \leq \inf T + \inf I + \inf F \leq \sup T + \sup I + \sup F \leq 3$, $\hspace{2cm}$ (107) we have a standard neutrosophic logic.

### 28. Definition of Extended Nonstandard Neutrosophic Logic

If in the above definition of general neutrosophic logic at least one of the neutrosophic components $T, I, F$ is a nonstandard real mobinad subset, neutrosophically included in or equal to the nonstandard real mobinad unit interval $]^-0, 1^+[_{MB}$, where $^-0 \leq_N \inf_N T + \inf_N I + \inf_N F \leq_N \sup_N T + \sup_N I + \sup_N F \leq 3^+$, $\hspace{1cm}$ (108) we have an extended nonstandard neutrosophic logic.

**Theorem 9**

If $M$ is a standard real set, $M \subset \mathbb{R}$, then $\inf_N(M) = \inf(M)$ and $\sup_N(M) = \sup(M)$. $\hspace{1cm}$ (109)

### 29. Definition of Extended General Neutrosophic Set

We extend and present in a clearer way our 1995 definition of neutrosophic set.

Let $\mathcal{U}$ be a universe of discourse of elements, and $S \in \mathcal{U}$ a subset.

A **Neutrosophic Set** is a set such that each element $x$ from $S$ has a degree of membership $(T)$, a degree of indeterminacy $(I)$, and a degree of nonmembership $(F)$, where $T, I, F$ are standard or nonstandard real mobinad subsets, neutrosophically included in or equal to the nonstandard real mobinat unit interval $]^-0, 1^+[_{MB}$,

with $T, I, F \subseteq_N ]^-0, 1^+[_{MB}$, (110)

where $^-0 \leq_N \inf_N T + \inf_N I + \inf_N F \leq_N \sup_N T + \sup_N I + \sup_N F \leq 3^+$. (111)

### 30. Definition of Standard Neutrosophic Set

If in the above general definition of neutrosophic set all neutrosophic components, $T, I, F$, are standard real subsets included in or equal to the classical real unit interval,

$T, I, F \subseteq [0, 1]$, where $0 \leq \inf T + \inf I + \inf F \leq \sup T + \sup I + \sup F \leq 3$, (112)

we have a standard neutrosophic set.

### 31. Definition of Extended Nonstandard Neutrosophic Set

If in the above general definition of neutrosophic set at least one of the neutrosophic components $T, I, F$ is a nonstandard real mobinad subsets, neutrosophically included in or equal to $]^-0, 1^+[_{MB}$,

where $^-0 \leq_N \inf_N T + \inf_N I + \inf_N F \leq_N \sup_N T + \sup_N I + \sup_N F \leq 3^+$, (113)

we have a nonstandard neutrosophic set.

### 32. Definition of Extended General Neutrosophic Probability

We extend and present in a clearer way our 1995 definition of neutrosophic probability.

Let $\mathcal{U}$ be a universe of discourse of events, and $E \in \mathcal{U}$ be an event.

A **Neutrosophic Probability** is a multivalued probability such that each event $E$ has a chance of occuring $(T)$, an indeterminate (unclear) chance of occuring or not occuring $(I)$, and a chance of not occuring $(F)$, where $T, I, F$ are standard or nonstandard real mobinad subsets, neutrosophically included in or equal to the nonstandard real mobinat unit interval $]^-0, 1^+[_{MB}$, $T, I, F \subseteq_N ]^-0, 1^+[_{MB}$, where $^-0 \leq_N \inf_N T + \inf_N I + \inf_N F \leq_N \sup_N T + \sup_N I + \sup_N F \leq 3^+$.

(114)

### 33. Definition of Standard Neutrosophic Probability

If in the above general definition of neutrosophic probability all neutrosophic components, $T, I, F$, are standard real subsets, included in or equal to the standard unit interval,

$T, I, F \subseteq [0, 1]$, where $0 \leq \inf T + \inf I + \inf F \leq \sup T + \sup I + \sup F \leq 3$, (115)

we have a standard neutrosophic probability.

### 34. Definition of Extended Nonstandard Neutrosophic Probability

If in the above general definition of neutrosophic probability at least one of the neutrosophic components $T$, $I$, $F$ is a nonstandard real mobinad subsets, neutrosophically included in or equal to $]^-0, 1^+[_{MB}$,

where $^-0 \leq_N \inf_N T + \inf_N I + \inf_N F \leq_N \sup_N T + \sup_N I + \sup_N F \leq 3^+$, (116)

we have a nonstandard neutrosophic probability.

### 35. Classical Operations with Real Sets

Let $A, B \subseteq \mathbb{R}$ be two real subsets. Let ⊛ and ∗ denote any of the <u>real subset classical operations</u> and <u>real number classical operations</u> respectively: addition (+), subtraction (−), multiplication (×), division (÷), and power (^).

Then, $A \circledast B = \{a * b, \text{where } a \in A \text{ and } b \in B\}$. (117)

Thus:

$$A \oplus B = \{a + b \mid a \in A, b \in B\}$$
$$A \ominus B = \{a - b \mid a \in A, b \in B\}$$
$$A \otimes B = \{a \times b \mid a \in A, b \in B\}$$
$$A \oslash B = \{a \div b \mid a \in A, b \in B - \{0\}\}$$
$$A^B = \{a \wedge b \mid a \in A, a > 0; b \in B\}$$

(118-122)

For the division (÷), of course, we consider $b \neq 0$. While for the power (^), we consider $a > 0$.

### 36. Operations on the Nonstandard Real MoBiNad Set ($NR_{MB}$)

For all nonstandard (addition, subtraction, multiplication, division, and power) operations,

for $\alpha, \beta \in_N NR_{MB}$, $\alpha *_N \beta =_N \mu_N(\alpha) \circledast \mu_N(\beta)$ (123)

where $*_N$ is any neutrosophic arithmetic operations with neutrosophic numbers (+N, -N, $\times_N, \div_N$, ^N), while the corresponding ⊛ is an arithmetic operation with real subsets.

So, we approximate the nonstandard operations by standard operations of real subsets.

We sink the nonstandard neutrosophic real mobinad operations into the standard real subset operations, then we resurface the last ones back to the nonstandard neutrosophic real mobinad set.

Let $\varepsilon_1$ and $\varepsilon_2$ be two non-null positive infinitesimals. We present below some particular cases, all others should be deduced analogously.

## Nonstandard Addition

*First Method*

$(^-a) + (^-b) =_N (a - \varepsilon_1, a) + (b - \varepsilon_2, b) =_N (a + b - \varepsilon_1 - \varepsilon_2, a + b) =_N (a + b - \varepsilon, a + b) =_N {}^-(a + b)$, where we denoted $\varepsilon_1 + \varepsilon_2 = \varepsilon$ (the addition of two infinitesimals is also an infinitesimal). (124)

*Second Method*

$(^-a) + (^-b) =_N (a - \varepsilon_1) + (b - \varepsilon_2) =_N (a + b - \varepsilon_1 - \varepsilon_2) =_N {}^-(a + b).$ (125)

Adding two left monads, one also gets a left monad.

## Nonstandard Subtraction

*First Method*

$(^-a) - (^-b) =_N (a - \varepsilon_1, a) - (b - \varepsilon_2, b) =_N (a - \varepsilon_1 - b, a - b + \varepsilon_2) =_N (a - b - \varepsilon_1, a - b + \varepsilon_2) =_N \begin{pmatrix} - & 0 & + \\ & a - b & \end{pmatrix}$ (126)

*Second Method*

$(^-a) - (^-b) =_N (a - \varepsilon_1) - (b - \varepsilon_2) =_N a - b - \varepsilon_1 + \varepsilon_2,$ (127)

since $\varepsilon_1$ and $\varepsilon_2$ may be any positive infinitesimals,

$$=_N \begin{cases} {}^-(a - b), \text{when } \varepsilon_1 > \varepsilon_2; \\ \begin{pmatrix} 0 \\ a - b \end{pmatrix}, \text{when } \varepsilon_1 = \varepsilon_2 \\ (a - b)^+, \text{when } \varepsilon_1 < \varepsilon_2. \end{cases} =_N \begin{pmatrix} 0 \\ a - b \end{pmatrix} =_N a - b;$$ (128-130)

Subtracting two left monads, one obtains an unpierced binad (that's why the unpierced binad had to be introduced).

## Nonstandard Division

Let $a, b > 0$.

$(^-a) \div (^-b) =_N (a - \varepsilon_1, a) \div (b - \varepsilon_2, b) =_N \left( \frac{a - \varepsilon_1}{b}, \frac{a}{b - \varepsilon_2} \right).$ (131)

Since $\varepsilon_1 > 0$ and $\varepsilon_2 > 0$, $\frac{a - \varepsilon_1}{b} < \frac{a}{b}$ and $\frac{a}{b - \varepsilon_2} > \frac{a}{b}$, (132)

while between $\frac{a - \varepsilon_1}{b}$ and $\frac{a}{b - \varepsilon_2}$ there is a continuum whence there are some infinitesimals $\varepsilon_1^0$ and $\varepsilon_2^0$ such that $\frac{a - \varepsilon_1^0}{b - \varepsilon_2^0} = \frac{a}{b}$, or $ab - b\varepsilon_1^0 = ab - a\varepsilon_2^0$, and for a given $\varepsilon_1^0$ there exist an $\varepsilon_2^0 = \varepsilon_1^0 \cdot \frac{b}{a}$. (133)

Whence $\frac{(^-a)}{(^-b)} =_N \left(^-\ 0\ ^+\atop \frac{a}{b}\right).$ (134)

For $a$ or/and $b$ negative numbers, it's similar but it's needed to compute the $inf_N$ and $sup_N$ of the products of intervals.

Dividing two left monads, one obtains an unpierced binad.

**Nonstandard Multiplication**

Let $a, b \geq 0$.

$(^-a^0) \times (^-b^{0\ +}) =_N (a - \varepsilon_1, a] \times (b - \varepsilon_2, b + \varepsilon_2) =_N ((a - \varepsilon_1) \cdot (b - \varepsilon_2), a \cdot (b + \varepsilon_2)) =_N (^-ab^{0\ +})$ (135)

since $(a - \varepsilon_1) \cdot (b - \varepsilon_2) < a \cdot b$ and $a \cdot (b + \varepsilon_2) > a \cdot b$. (136)

For $a$ or/and $b$ negative numbers, it's similar but it's needed to compute the $inf_N$ and $sup_N$ of the products of intervals.

Multiplying a positive left monad closed to the right, with a positive unpierced binad, one obtains an unpierced binad.

**Nonstandard Power**

Let $a, b > 1$.

$(^0a^+)^{(^-b^0)} =_N [a, a + \varepsilon_1)^{(b - \varepsilon_2, b]} =_N (a^{b-\varepsilon_2}, (a + \varepsilon_1)^b) =_N \left(^-\ 0\ ^+\atop a^b\right)$ (137)

since $a^{b-\varepsilon_1} < a^b$ and $(a + \varepsilon_1)^b > a^b$. (138)

Raising a right monad closed to the left to a power equal to a left monad closed to the right, for both monads above 1, the result is an unpierced binad.

### Consequence

In general, when doing arithmetic operations on nonstandard real monads and binads, the result may be a different type of monad or binad.

That's why is was imperious to extend the monads to closed monads, and the pierced binad to unpierced binad, in order to have the whole nonstandard neutrosophic real mobinad set closed under arithmetic operations.

## 37. Conditions of Neutrosophic Nonstandard Inequalities

Let $NR_{MB}$ be the Nonstandard Real MoBiNad. Let's endow $(NR_{MB}, <_N)$ with a neutrosophic inequality.

Let $\alpha, \beta \in NR_{MB}$, where $\alpha, \beta$ may be real numbers, monads, or binads.

And let

$$\left(\overset{-}{a}\right), \left(\overset{-0}{a}\right), \left(\overset{+}{a}\right), \left(\overset{0+}{a}\right), \left(\overset{-+}{a}\right), \left(\overset{-0+}{a}\right) \in NR_{MB}, \text{ and } \left(\overset{-}{b}\right), \left(\overset{-0}{b}\right), \left(\overset{+}{b}\right), \left(\overset{0+}{b}\right), \left(\overset{-+}{b}\right), \left(\overset{-0+}{b}\right) \in NR_{MB}, \quad (139)$$

be the left monads, left monads closed to the right, right monads, right monads closed to the left, and binads, and binads nor prierced of the elements (standard real numbers) *a* and *b* respectively. Since all monads and binads are real subsets, we may treat the single real numbers

$a = [a, a]$ and $b = [b, b]$ as real subsets too. $\quad (140)$

$NR_{MB}$ is a set of subsets, and thus we deal with neutrosophic inequalities between subsets.

  i) If the subset $\alpha$ has many of its elements above all elements of the subset $\beta$, then $\alpha >_N \beta$ (partially).
  ii) If the subset $\alpha$ has many of its elements below all elements of the subset $\beta$, then $\alpha <_N \beta$ (partially).
  iii) If the subset $\alpha$ has many of its elements equal with elements of the subset $\beta$, then $\alpha =_N \beta$ (partially).

If the subset $\alpha$ verifies *i)* and *iii)* with respect to subset $\beta$, then $\alpha \geq_N \beta$.

If the subset $\alpha$ verifies *ii)* and *iii)* with respect to subset $\beta$, then $\alpha \leq_N \beta$.

If the subset $\alpha$ verifies *i)* and *ii)* with respect to subset $\beta$, then there is no neutrosophic order (inequality) between $\alpha$ and $\beta$.

{ For example, between *a* and *($^-a^+$)* there is no neutrosophic order, similarly between *a* and $\overset{-0+}{a}$. }

Similarly, if the subset $\alpha$ verifies *i), ii)* and *iii)* with respect to subset $\beta$, then there is no neutrosophic order (inequality) between $\alpha$ and $\beta$.

## 38. Open Neutrosophic Research

The quantity or measure of "many of its elements" of the above *i), ii),* or *iii)* conditions depends on each *neutrosophic application* and on its *neutrosophic experts*.

An approach would be to employ the *Neutrosophic Measure* [21, 22], that handles indeterminacy, which may be adjusted and used in these cases.

In general, we do not try in purpose to validate or invalidate an existing scientific result, but to investigate how an existing scientific result behaves in a new environment (that may contain indeterminacy), or in a new application, or in a new interpretation.

### 39. Nonstandard Neutrosophic Inequalities

For the *neutrosophic nonstandard inequalities*, we propose based on the previous six neutrosophic equalities, the following:

$$(^-a) <_N a <_N (a^+) \qquad (141)$$

since the standard real interval $(a - \varepsilon, a)$ is below $a$, and $a$ is below the standard real interval $(a, a + \varepsilon)$ by using the approximation provided by the nonstandard neutrosophic function $\mu$,

or because $\forall x \in R_+^*, a - x < a < a + x,$ \qquad (142)

where $x$ is of course a (nonzero) positive infinitesimal (the above double neutrosophic inequality actually becomes a double classical standard real inequality for each fixed positive infinitesimal).

The converse double neutrosophic inequality is also neutrosophically true:

$$(a^+) >_N a >_N (^-a) \qquad (143)$$

Another nonstandard neutrosophic double inequality:

$$(^-a) \leq_N (^-a^+) \leq_N (a^+) \qquad (144)$$

This double neutrosophic inequality may be justified since $(^-a^+) = (^-a) \cup (a^+)$ and, geometrically, on the Real Number Line, the number $a$ is in between the subsets $^-a = (a-\varepsilon, a)$ and $a^+ = (a, a+\varepsilon)$, so:

$$(^-a) \leq_N (^-a) \cup (a^+) \leq_N (a^+) \qquad (145)$$

whence the left side of the inequality's middle term coincides with the inequality first term, while the right side of the inequality middle term coincides with the third inequality term.

Conversely, it is neutrosophically true as well:

$$(a^+) \geq_N (^-a) \cup (a^+) \geq_N (^-a) \qquad (146)$$

Also, $^-a \leq_N {}^{-0}a \leq_N a \leq_N {}^{0+}a \leq_N {}^+a$ and $^-a \leq_N {}^{-+}a \leq_N {}^{-0+}a \leq_N {}^+a$. \qquad (147)

Conversely, they are also neutrosophically true:

$${}^+a \geq_N {}^{0+}a \geq_N a \geq_N {}^{-0}a \geq_N {}^-a \text{ and } {}^+a \geq_N {}^{-0+}a \geq_N {}^{-+}a \geq_N {}^-a \text{ respectively.} \qquad (148)$$

If $a > b$, which is a (standard) classical real inequality, then we have the following neutrosophic nonstandard inequalities:

$$a >_N (^-b), \quad a >_N (b^+), \quad a >_N (^-b^+), \quad a >_N {}^{-0}b, a >_N {}^{0+}b, a >_N {}^{-0+}b; \qquad (149)$$

$(^-a) >_N b$, $(^-a) >_N (^-b)$, $(^-a) >_N (b^+)$, $(^-a) >_N (^-b^+)$, $\overset{-}{a} >_N \overset{-0-}{b}, \overset{}{a} >_N \overset{0+-}{b}, \overset{}{a} >_N \overset{-0+}{b}$; (150)

$(a^+) >_N b$, $(a^+) >_N (^-b)$, $(a^+) >_N (b^+)$, $(a^+) >_N (^-b^+)$, $\overset{+}{a} >_N \overset{-0+}{b}, \overset{}{a} >_N \overset{0++}{b}, \overset{}{a} >_N \overset{-0+}{b}$; (151)

$(^-a^+) >_N b$, $(^-a^+) >_N (^-b)$, $(^-a^+) >_N (b^+)$, $(^-a^+) >_N (^-b^+)$, etc. (152)

**No Ordering Relationships**

For any standard real number $a$, there is no relationship of order between

the elements $a$ and $(^-a^+)$, nor between the elements $a$ and $\left(\overset{-0+}{a}\right)$. (153)

Therefore, $NR_{MB}$ is a neutrosophically partially order set.

If one removes all binads from $NR_{MB}$, then $(NR_{MB}, \leq_N)$ is neutrosophically totally ordered. (154)

**Theorem 10.**

Using the nonstandard general notation one has:

If $a > b$, which is a (standard) classical real inequality, then

$\overset{m_1}{a} >_N \overset{m_2}{b}$ for any $m_1, m_2 \in \{\ ,\ ^-,\ ^{-0},\ ^+,\ ^{+0},\ ^{-+},\ ^{-0+}\}$. (155)

And conversely,

If $a < b$, which is a (standard) classical real inequality, then

$\overset{m_1}{a} <_N \overset{m_2}{b}$ for any $m_1, m_2 \in \{\ ,\ ^-,\ ^{-0},\ ^+,\ ^{+0},\ ^{-+},\ ^{-0+}\}$. (156)

### 40. Nonstandard Neutrosophic Equalities

Let $a, b$ be standard real numbers; if $a = b$ that is a (classical) standard equality, then:

$(^-a) =_N (^-b)$, $(a^+) =_N (b^+)$, $(^-a^+) =_N (^-b^+)$, (157)

$\left(\overset{-0}{a}\right) =_N \left(\overset{-0}{b}\right), \left(\overset{0+}{a}\right) =_N \left(\overset{0+}{b}\right), \left(\overset{-0+}{a}\right) =_N \left(\overset{-0+}{b}\right)$. (158)

### 41. Nonstandard Neutrosophic Belongingness

On the nonstandard real set $NR_{MB}$, we say that

$\overset{m}{c} \in_N ]\overset{m_1}{a}, \overset{m_2}{b}[$ iff $\overset{m_1}{a} \leq_N \overset{m}{c} \leq_N \overset{m_2}{b}$, (159)

where $m_1, m_2, m \in \{\ ,\ ^-,\ ^{-0},\ ^+,\ ^{+0},\ ^{-+},\ ^{-0+}\}$. (160)

{ We use the previous nonstandard neutrosophic inequalities. }

### 42. Nonstandard Hesitant Sets

Are sets of the form: $A = \{a_1, a_2, ..., a_n\}$, $2 \leq n < \infty$, $A \subset_N NR_{MB}$,

where at least one element $a_{i_0}, 1 \leq i_0 \leq n$, is either an infinitesimal, or a monad or a binad (of any type);

while other elements may be standard real numbers, infinitesimals, or also monads or binads (of any type). (161)

If the neutrosophic components T, I, F are nonstandard hesitant sets, then one has a

*Nonstandard Hesitant Neutrosophic Logic / Set / Probability*.

### 43. Nonstandard Neutrosophic Strict Interval Inclusion

On the nonstandard real set $NR_{MB}$,

$] {}^{m_1}a, {}^{m_2}b [ \subset_N ] {}^{m_3}c, {}^{m_4}d [$ iff (162)

${}^{m_3}c \leq_N {}^{m_1}a <_N {}^{m_2}b <_N {}^{m_4}d$ or ${}^{m_3}c <_N {}^{m_1}a <_N {}^{m_2}b \leq_N {}^{m_4}d$ or ${}^{m_3}c <_N {}^{m_1}a <_N {}^{m_2}b <_N {}^{m_4}d$. (163)

### 44. Nonstandard Neutrosophic (Non-Strict) Interval Inclusion

On the nonstandard real set $NR_{MB}$,

$] {}^{m_1}a, {}^{m_2}b [ \subseteq_N ] {}^{m_3}c, {}^{m_4}d [$ iff (164)

${}^{m_3}c \leq_N {}^{m_1}a <_N {}^{m_2}b \leq_N {}^{m_4}d$. (165)

### 45. Nonstandard Neutrosophic Strict Set Inclusion

The nonstandard set $A$ is neutrosophically strictly included in the nonstandard set $B$, $A \subset_N B$, if:

$\forall x \in_N A, x \in_N B$, and $\exists y \in_N B : y \notin_N A$. (166)

### 46. Nonstandard Neutrosophic (Non-Strict) Set Inclusion

The nonstandard set $A$ is neutrosophically not-strictly included in the nonstandard set $B$,

$A \subseteq_N B$, iff: (167)

$\forall x \in_N A, x \in_N B$. (168)

## 47. Nonstandard Neutrosophic Set Equality

The nonstandard sets $A$ and $B$ are neutrosophically equal, $A =_N B$, iff: (169)

$A \subseteq_N B$ and $B \subseteq_N A$. (170)

## 48. The Fuzzy, Neutrosophic, and Plithogenic Logical Connectives ∧, ∨, →

All fuzzy, intuitionistic fuzzy, and neutrosophic logic operators are *inferential approximations*, not written in stone. They are improved from application to application.

Let's denote:

$\wedge_F, \wedge_N, \wedge_P$ representing respectively the fuzzy conjunction, neutrosophic conjunction, and plithogenic conjunction; (171)

similarly

$\vee_F, \vee_N, \vee_P$ representing respectively the fuzzy disjunction, neutrosophic disjunction, and plithogenic disjunction, (172)

and

$\rightarrow_F, \rightarrow_N, \rightarrow_P$ representing respectively the fuzzy implication, neutrosophic implication, and plithogenic implication. (173)

I agree that my beginning neutrosophic operators (when I applied the same *fuzzy t-norm*, or the same *fuzzy t-conorm*, to all neutrosophic components *T, I, F*) were less accurate than others developed later by the neutrosophic community researchers. This was pointed out since 2002 by Ashbacher [9] and confirmed in 2008 by Rivieccio [10]. They observed that if on $T_1$ and $T_2$ one applies a *fuzzy t-norm*, on their opposites $F_1$ and $F_2$ one needs to apply the *fuzzy t-conorm* (the opposite of fuzzy t-norm), and reciprocally.

About inferring $I_1$ and $I_2$, some researchers combined them in the same directions as $T_1$ and $T_2$.

Then:

$(T_1, I_1, F_1) \wedge_N (T_2, I_2, F_2) = (T_1 \wedge_F T_2, I_1 \wedge_F I_2, F_1 \vee_F F_2)$, (174)

$(T_1, I_1, F_1) \vee_N (T_2, I_2, F_2) = (T_1 \vee_F T_2, I_1 \vee_F I_2, F_1 \wedge_F F_2)$, (175)

$(T_1, I_1, F_1) \rightarrow_N (T_2, I_2, F_2) = (F_1, I_1, T_1) \vee_N (T_2, I_2, F_2) = (F_1 \vee_F T_2, I_1 \vee_F I_2, T_1 \wedge_F F_2);$ (176)

others combined $I_1$ and $I_2$ in the same direction as $F_1$ and $F_2$ (since both *I* and *F* are negatively qualitative neutrosophic components, while F is qualitatively positive neutrosophic component), the most used one:

$$(T_1, I_1, F_1) \wedge_N (T_2, I_2, F_2) = (T_1 \wedge_F T_2, I_1 \vee_F I_2, F_1 \vee_F F_2), \tag{177}$$

$$(T_1, I_1, F_1) \vee_N (T_2, I_2, F_2) = (T_1 \vee_F T_2, I_1 \wedge_F I_2, F_1 \wedge_F F_2), \tag{178}$$

$$(T_1, I_1, F_1) \to_N (T_2, I_2, F_2) = (F_1, I_1, T_1) \vee_N (T_2, I_2, F_2) = (F_1 \vee_F T_2, I_1 \wedge_F I_2, T_1 \wedge_F F_2). \tag{179}$$

Even more, recently, in an extension of neutrosophic set to *plithogenic set* [11] (which is a set whose each element is characterized by many attribute values), the *degrees of contradiction* $c(\,,\,)$ between the neutrosophic components T, I, F have been defined (in order to facilitate the design of the aggregation operators), as follows:

$c(T, F) = 1$ *(or 100%, because they are totally opposite)*, $c(T, I) = c(F, I) = 0.5$ *(or 50%, because they are only half opposite)*, (180)
then:

$$(T_1, I_1, F_1) \wedge_P (T_2, I_2, F_2) = (T_1 \wedge_F T_2, 0.5(I_1 \wedge_F I_2) + 0.5(I_1 \vee_F I_2), F_1 \vee_F F_2), \tag{181}$$

$$(T_1, I_1, F_1) \vee_P (T_2, I_2, F_2) = (T_1 \vee_F T_2, 0.5(I_1 \vee_F I_2) + 0.5(I_1 \wedge_F I_2), F_1 \wedge_F F_2). \tag{182}$$

$$(T_1, I_1, F_1) \to_N (T_2, I_2, F_2) = (F_1, I_1, T_1) \vee_N (T_2, I_2, F_2)$$
$$= (F_1 \vee_F T_2, 0.5(I_1 \vee_F I_2) + 0.5(I_1 \wedge_F I_2), T_1 \wedge_F F_2). \tag{183}$$

### 49. Fuzzy t-norms and Fuzzy t-conorms
The most used $\wedge_F$ (Fuzzy t-norms), and $\vee_F$ (Fuzzy t-conorms) are:

Let $a, b \in [0, 1]$. (184)

Fuzzy t-norms (fuzzy conjunctions, or fuzzy intersections):
$a \wedge_F b = min\{a, b\}$; (185)
$a \wedge_F b = ab$; (186)
$a \wedge_F b = max\{a + b - 1, 0\}$. (187)

Fuzzy t-conorms (fuzzy disjunctions, or fuzzy unions):
$a \vee_F b = max\{a, b\}$; (188)
$a \vee_F b = a + b - ab$; (189)
$a \vee_F b = min\{a + b, 1\}$. (190)

### 50. Nonstandard Neutrosophic Operators

**Nonstandard Neutrosophic Conjunctions**

$(T_1, I_1, F_1) \wedge_N (T_2, I_2, F_2) = (T_1 \wedge_F T_2, I_1 \vee_F I_2, F_1 \vee_F F_2) =$
$(\,inf_N(T_1, T_2),\, sup_N(I_1, I_2),\, sup_N(F_1, F_2)\,)$ (191)

$(T_1, I_1, F_1) \wedge_N (T_2, I_2, F_2) = (T_1 \wedge_F T_2, I_1 \vee_F I_2, F_1 \vee_F F_2) =$

$$( T_1 \times_N T_2, I_1 +_N I_2 -_N I_1 \times_N I_2, F_1 +_N F_2 -_N F_1 \times_N F_2 ) \tag{192}$$

**Nonstandard Neutrosophic Disjunctions**

$$(T_1, I_1, F_1) \vee_N (T_2, I_2, F_2) = (T_1 \vee_F T_2, I_1 \wedge_F I_2, F_1 \wedge_F F_2) =$$
$$( \sup_N(T_1, T_2), \inf_N(I_1, I_2), \inf_N(F_1, F_2) ) \tag{193}$$

$$(T_1, I_1, F_1) \vee_N (T_2, I_2, F_2) = (T_1 \vee_F T_2, I_1 \wedge_F I_2, F_1 \wedge_F F_2) =$$
$$( T_1 +_N T_2 -_N T_1 \times_N T_2, I_1 \times_N I_2, F_1 \times_N F_2 ) \tag{194}$$

**Nonstandard Neutrosophic Negations**

$$\neg(T_1, I_1, F_1) = (F_1, I_1, T_1) \tag{195}$$

$$\neg(T_1, I_1, F_1) = ( F_1, (1^+) -_N I_1, T_1 ) \tag{196}$$

**Nonstandard Neutrosophic Implications**

$$(T_1, I_1, F_1) \rightarrow_N (T_2, I_2, F_2) = (F_1, I_1, T_1) \vee_N (T_2, I_2, F_2) = (F_1 \vee_F T_2, I_1 \wedge_F I_2, T_1 \wedge_F F_2)$$
$$= ( F_1 +_N T_2 -_N F_1 \times_N T_2, I_1 \times_N I_2, T_1 \times_N F_2 ) \tag{197}$$

$$(T_1, I_1, F_1) \rightarrow_N (T_2, I_2, F_2) = (F_1, (1^+) -_N I_1, T_1) \vee_N (T_2, I_2, F_2)$$
$$= (F_1 \vee_F T_2, ((1^+) -_N I_1) \wedge_F I_2, T_1 \wedge_F F_2) = ( F_1 +_N T_2 -_N F_1 \times_N T_2, ((1^+) -_N I_1) \times_N I_2, T_1 \times_N F_2 ) \tag{198}$$

Let $P_1(T_1, I_1, F_1)$ and $P_2(T_2, I_2, F_2)$ be two nonstandard neutrosophic logical propositions, whose nonstandard neutrosophic components are respectively: $T_1, I_1, F_1, T_2, I_2, F_2 \in_N NR_{MB}$. (199)

## *51. Numerical Examples of Nonstandard Neutrosophic Operators*

Let's take a particular numeric example, where:

$$P_1 =_N (\overset{0+}{0.3}, \overset{-+}{0.2}, 0.4), P_2 =_N (\overset{-0}{0.6}, \overset{-0+}{0.1}, \overset{+}{0.5}) \tag{200}$$

are two nonstandard neutrosophic logical propositions.
We use the nonstandard arithmetic operations previously defined

*Numerical Example of Nonstandard Neutrosophic Conjunction*

$$\overset{0+}{0.3} \times_N \overset{-0}{0.6} =_N [0.3, 0.3+\varepsilon_1) \times (0.6-\varepsilon_2, 0.6] = (0.18-0.3\varepsilon_2, 0.18+0.6\varepsilon_1) =_N \overset{-0+}{0.18} \tag{201}$$

$$0.2+_N^{-+} \ 0.1-_N^{-0+} \ 0.2\times_N^{-+} \ 0.1 =_N [(0.2-\varepsilon_1,0.2)\cup(0.2,0.2+\varepsilon_1)]+(0.1-\varepsilon_2,0.1+\varepsilon_2)$$
$$-[(0.2-\varepsilon_1,0.2)\cup(0.2,0.2+\varepsilon_1)]\times(0.1-\varepsilon_2,0.1+\varepsilon_2)$$
$$=[(0.3-\varepsilon_1-\varepsilon_2,0.3+\varepsilon_2)\cup(0.3-\varepsilon_2,0.3+\varepsilon_1+\varepsilon_2)]$$
$$-[(0.2-\varepsilon_1)\times(0.1-\varepsilon_2),(0.02+0.2\varepsilon_2)]\cup[(0.02-0.2\varepsilon_2),(0.2+\varepsilon_1)\times(0.1+\varepsilon_2)]$$
$$=[\overset{-0+}{0.3}\cup\overset{-0+}{0.3}]-[\overset{-0+}{0.02}\cup\overset{-0+}{0.02}]=[\overset{-0+}{0.3}]-[\overset{-0+}{0.02}]=\overset{-0+}{0.3}-\overset{-0+}{0.02}=_N \overset{-0+}{0.28} \quad (202)$$

$$0.4+_N \overset{+}{0.5}=_N [0.4,0.4]+(0.5,0.5+\varepsilon_1)-[0.4,0.4]\times(0.5,0.5+\varepsilon_1)$$
$$=(0.4+0.5,0.4+0.5+\varepsilon_1)-(0.4\times0.5,0.4\times0.5+0.4\varepsilon_1)$$
$$=(0.9,0.9+\varepsilon_1)-(0.2,0.2+0.4\varepsilon_1)$$
$$=(0.9-0.2-0.4\varepsilon_1,0.9+\varepsilon_1-0.2)=(0.7-0.4\varepsilon_1,0.7+\varepsilon_1)=_N \overset{-0+}{0.70} \quad (203)$$

Whence

$$P_1 \wedge_N P_2 =_N (\overset{-0+}{0.18},\overset{-0+}{0.28},\overset{-0+}{0.70}) \quad (204)$$

*Numerical Example of Nonstandard Neutrosophic Disjunction*

$$\overset{0+}{0.3}+_N \overset{-0}{0.6}-\overset{0+}{0.3}\times_N \overset{-0}{0.6}=_N \{[0.3,0.3+\varepsilon_1)+(0.6-\varepsilon_1,0.6]\}-\{[0.3,0.3+\varepsilon_1)\times(0.6-\varepsilon_1,0.6]\}$$
$$=(0.9-\varepsilon_1,0.9+\varepsilon_1)-(0.18-0.3\varepsilon_1,0.18+0.6\varepsilon_1)=(0.72-1.6\varepsilon_1,0.72+1.3\varepsilon_1)=_N \overset{-0+}{0.72} \quad (205)$$

$$\overset{-+}{0.2}\times_N \overset{-0+}{0.1}=_N \left(\overset{-0+}{0.2\times0.1}\right)=_N \overset{-0+}{0.02} \quad (206)$$

$$0.4\times_N \overset{+}{0.5}=_N \left(\overset{+}{0.4\times0.5}\right)=_N \overset{+}{0.20} \quad (207)$$

Whence $P_1 \vee_N P_2 =_N (\overset{-0+}{0.72},\overset{-0+}{0.28},\overset{+}{0.20}) \quad (208)$

*Numerical Example of Nonstandard Neutrosophic Negation*

$$\neg_N P_1 =_N \neg_N(\overset{0+}{0.3},\overset{-+}{0.2},0.4)=_N (0.4,\overset{-+}{0.2},\overset{0+}{0.3}) \quad (209)$$

*Numerical Example of Nonstandard Neutrosophic Implication*

$$(P_1 \rightarrow_N P_2) \Leftrightarrow_N (\neg_N P_1 \vee_N P_2) =_N (0.4,\overset{-+}{0.2},\overset{0+}{0.3}) \vee_N (\overset{-0}{0.6},\overset{-0+}{0.1},\overset{+}{0.5}) \quad (210)$$

Afterwards,

$$0.4 +_N \overset{-0}{0.6} -_N 0.4 \times_N \overset{-0}{0.6} =_N \left(0.4 \overset{-0}{+} 0.6\right) -_N \left(0.4 \overset{-0}{\times} 0.6\right) =_N \overset{-0}{1.0} -_N \overset{-0}{0.24} =_N \overset{-0+}{0.76}$$

$$\overset{-+}{0.2} \times_N \overset{-0+}{0.1} =_N \overset{-0+}{0.02}$$

$$\overset{0+}{0.3} \times_N \overset{+}{0.5} =_N \overset{+}{0.15}$$

(211-213)

whence

$$\neg_N P_1 =_N (\overset{-0+}{0.76}, \overset{-0+}{0.02}, \overset{+}{0.15}).$$

(214)

### Conclusion

In the history of mathematics, critics on nonstandard analysis, in general, have been made by Paul Halmos, Errett Bishop, Alain Connes and others.

That's why we have extended in 1998 for the first time the monads to *pierced binad*, and then in 2019 for the second time we extended the left monad to *left monad closed to the right*, the right monad to *right monad closed to the left*, and the pierced binad to *unpierced binad*. These were necessary in order to construct a general nonstandard neutrosophic real mobinad space, which is closed under the nonstandard neutrosophic arithmetic operations (such as addition, subtraction, multiplication, division, and power) which are needed in order to be able to define the nonstandard neutrosophic operators (such as conjunction, disjunction, negation, implication, equivalence) on this space, and to transform the newly constructed nonstandard neutrosophic real mobinad space into a lattice of first order (as partially ordered nonstandard set, under the neutrosophic inequality $\leq_N$) and a lattice of second type [as algebraic structure, endowed with two binary laws: neutrosophic infimum (*infN*) and neutrosophic supremum (*supN*)].

### References

[1] Takura Imamura, *Note on the Definition of Neutrosophic Logic*, arxiv.org, 7 Nov. 2018.
[2] Xindong Peng and Jingguo Dai, *A bibliometric analysis of neutrosophic set: two decades review from 1998 to 2017*, Artificial Intelligence Review, Springer, 18 August 2018; http://fs.unm.edu/BibliometricNeutrosophy.pdf
[3] Florentin Smarandache, *n-Valued Refined Neutrosophic Logic and Its Applications in Physics*, Progress in Physics, 143-146, Vol. 4, 2013;
http://fs.unm.edu/n-ValuedNeutrosophicLogic-PiP.pdf
[4] F. Smarandache, *Neutrosophy, A New Branch of Philosophy*, <Multiple Valued Logic / An International Journal>, USA, ISSN 1023-6627, Vol. 8, No. 3, pp. 297-384, 2002.
[5] Florentin Smarandache, *Neutrosophy. / Neutrosophic Probability, Set, and Logic*, ProQuest Information & Learning, Ann Arbor, Michigan, USA, 105 p., 1998;
http://fs.unm.edu/eBook-Neutroosphics6.pdf.

[6] F. Smarandache, *A Unifying Field in Logics: Neutrosophic Logic*, <Multiple Valued Logic / An International Journal>, USA, ISSN 1023-6627, Vol. 8, No. 3, pp. 385-438, 2002. {The whole issue of this journal is dedicated to Neutrosophy and Neutrosophic Logic.}
[7] Florentin Smarandache, *Definition of neutrosophic logic — a generalization of the intuitionistic fuzzy logic*, Proceedings of the 3rd Conference of the European Society for Fuzzy Logic and Technology, 2003, pp. 141–146.
[8] Florentin Smarandache, Neutrosophic Overset, Neutrosophic Underset, and Neutrosophic Offset. Similarly for Neutrosophic Over-/Under-/Off- Logic, Probability, and Statistics, 168 p., Pons Editions, Bruxelles, Belgique, 2016; https://arxiv.org/ftp/arxiv/papers/1607/1607.00234.pdf
[9] Charles Ashbacher, Introduction to Neutrosophic Logic, ProQuest Information & Learning, Ann Arbor, 2002, http://fs.unm.edu/IntrodNeutLogic.pdf
[10] Umberto Rivieccio, Neutrosophic logics: Prospects and problems, Fuzzy Sets and Systems, v. 159, issue 14, 1860–1868, 2008.
[11] Plithogeny, Plithogenic Set, Logic, Probability, and Statistics, by Florentin Smarandache, Pons Publishing House, Brussels, Belgium, 141 p., 2017;
arXiv.org (Cornell University), Computer Science - Artificial Intelligence, 03Bxx: https://arxiv.org/ftp/arxiv/papers/1808/1808.03948.pdf
[12] Nguyen Xuan Thao, Florentin Smarandache, *(I, T)-Standard neutrosophic rough set and its topologies properties*, Neutrosophic Sets and Systems, Vol. 14, 2016, pp. 65-70; doi.org/10.5281/zenodo.570892
http://fs.unm.edu/NSS/RoughStandardNeutrosophicSets.pdf
[13] Nguyen Xuan Thao, Bui Cong Cuong, Florentin Smarandache, *Rough Standard Neutrosophic Sets: An Application on Standard Neutrosophic Information Systems*, Neutrosophic Sets and Systems, Vol. 14, 2016, pp. 80-92; doi.org/10.5281/zenodo.570890
http://fs.unm.edu/NSS/RoughStandardNeutrosophicSets.pdf
[14] Bui Cong Cuong, Pham Hong Phong, Florentin Smarandache, *Standard Neutrosophic Soft Theory - Some First Results*, Neutrosophic Sets and Systems, Vol. 12, 2016, pp. 80-91; doi.org/10.5281/zenodo.571149
http://fs.unm.edu/NSS/StandardNeutrosophicSoftTheory.pdf
[15] Insall, Matt and Weisstein, Eric W. "Nonstandard Analysis." From *MathWorld*--A Wolfram Web Resource. http://mathworld.wolfram.com/NonstandardAnalysis.html
[16] Insall, Matt. "Transfer Principle." From *MathWorld*--A Wolfram Web Resource, created by Eric W. Weisstein. http://mathworld.wolfram.com/TransferPrinciple.html
[17] F. Smarandache, *Applications of Neutrosophic Sets in Image Identification, Medical Diagnosis, Fingerprints and Face Recognition* and *Neutrosophic Overset/Underset/Offset*, COMSATS Institute of Information Technology, Abbottabad, Pakistan, December 26th, 2017;
[18] F. Smarandache, *Interval-Valued Neutrosophic Oversets, Neutrosophic Understes, and Neutrosophic Offsets*, International Journal of Science and Engineering Investigations, Vol. 5, Issue 54, 1-4, July 2016.
[19] F. Smarandache, *Operators on Single-Valued Neutrosophic Oversets, Neutrosophic Undersets, and Neutrosophic Offsets*, Journal of Mathematics and Informatics, Vol. 5, 63-67, 2016.
[20] Florentin Smarandache, *Answers to Imamura Note on the Definition of Neutrosophic Logic*, 24 November 2018, http://arxiv.org/abs/1812.02534,
https://arxiv.org/ftp/arxiv/papers/1812/1812.02534.pdf
[21] Florentin Smarandache, *Introduction to Neutrosophic Measure, Neutrosophic Integral, and*

*Neutrosophic Probability*, Sitech & Educational, Craiova, Columbus, 140 p., 2013; https://arxiv.org/ftp/arxiv/papers/1311/1311.7139.pdf

[22] Florentin Smarandache, An Introduction to Neutrosophic Measure, American Physical Society, April Meeting 2014, Volume 59, Number 5, Saturday-Tuesday, April 5-8, 2014; Savannah, Georgia, http://meetings.aps.org/Meeting/APR14/Session/T1.26

[23] F. Smarandache, A*bout Nonstandard Neutrosophic Logic (Answers to Imamura's 'Note on the Definition of Neutrosophic Logic')*, Cornell University, New York City, USA, (Submitted on 24 Nov 2018 (v1), last revised 13 Feb 2019 (this version, v2)),
Abstract: https://arxiv.org/abs/1812.02534v2
Full Paper: https://arxiv.org/ftp/arxiv/papers/1812/1812.02534.pdf